\documentclass[a4paper,reqno,11pt]{amsart}

\usepackage[utf8]{inputenc}
\usepackage{microtype}

\usepackage[margin=1.1in]{geometry}

\usepackage{amssymb}
\usepackage{mathrsfs}
\usepackage{mathtools}
\usepackage{enumitem}
\usepackage{comment}

\usepackage{color}
\definecolor{gray}{rgb}{.75,.75,.75}
\usepackage{tikz}
\usetikzlibrary{topaths}
\usetikzlibrary{tikzmark}
\usetikzlibrary{intersections}
\usetikzlibrary{decorations.pathreplacing}

\usepackage{savesym}
\usepackage[all]{xy}
\savesymbol{cir}

\usepackage{hyperref}

\newtheorem{theorem}{Theorem}[section]

\theoremstyle{definition}
\newtheorem{definition}[theorem]{Definition}

\newtheorem{example}[theorem]{Example}

\newcommand{\Z}{\mathbb{Z}}
\newcommand{\N}{\mathbb{N}}
\newcommand{\Q}{\mathbb{Q}}

\newcommand{\calO}{\mathcal{O}}

\DeclareMathOperator{\id}{id}

\newcommand{\abels}{\mathit{Ab}}

\newcommand{\F}{\mathrm{F}}

\newcommand{\clone}{\kappa}
\newcommand{\symmclone}{\varsigma}
\newcommand{\colim}{\varinjlim}

\newcommand{\image}{\operatorname{im}}

\newcommand{\defeq}{\mathbin{\vcentcolon =}}

\usepackage{ifthen}
\usepackage{xargs}
\newcommand{\optionalarg}[2]{
\ifthenelse{\equal{#2}{}}{%
#1}{%
#1(#2)}
}

\newcommand{\Thomp}[1]%
   {\optionalarg{\mathscr{T}}{#1}}                 
   
\newcommand{\Thkern}[1]%
   {\optionalarg{\mathscr{K}}{#1}}                 

\newcommand{\Stein}[1]%
   {\optionalarg{\mathscr{X}}{#1}}                 

\newcommand{\Fbr}%
   {F_{\operatorname{br}}}                 

\newcommand{\Vbr}%
   {V_{\operatorname{br}}}                 
   
\newcommand{\Vmock}%
   {V_{\operatorname{mock}}}                 
   
\newcommand{\Vloop}%
   {V_{\operatorname{loop}}}                 
   
\newcommand{\Floop}%
   {F_{\operatorname{loop}}}                 

\begin{document}\parindent=0pt

\title{A user's guide to cloning systems}
\date{\today}


\author[M.~C.~B.~Zaremsky]{Matthew C.~B.~Zaremsky}
\address{Department of Mathematics, Cornell University, Ithaca, NY 14853}
\email{zaremskym@gmail.com}

\begin{abstract}
 In joint work of the author with Stefan Witzel, a procedure was developed for building new examples of groups in the extended family of R.~Thompson's groups, using what we termed \emph{cloning systems}. These new Thompson-like groups can be thought of as limits of families of groups, though unlike other limiting processes, e.g., direct limits, these tend to be well behaved with respect to finiteness properties. In this expository note, we distill the crucial parts of that 50-page paper into a more digestible form, for those curious to understand the construction but less curious about the gritty details. We also give some new examples, involving signed symmetric groups and twisted braid groups.
\end{abstract}

\maketitle
\thispagestyle{empty}

\section{Introduction}\label{sec:intro}

The notion of a \emph{cloning system} on a family of groups $(G_n)_{n\in\N}$ was introduced by Stefan Witzel and the author, in the paper \cite{witzel17}. Given a cloning system on $(G_n)_{n\in\N}$, one gets a group $\Thomp{G_*}$, called the \emph{generalized Thompson group} for the cloning system (more often called a \emph{Thompson-like group}). One original motivation for axiomatizing the cloning system construction was to build a general framework giving rise to various preexisting versions of the R.~Thompson groups, for example groups called $F$, $V$, $\Vbr$ and $\Fbr$ (using the families $(\{1\})$, $(S_n)$, $(B_n)$ and $(PB_n)$ respectively), and also some new examples the authors found, for example using the family $(B_n(R))$. Here $B_n(R)$ is the group of upper triangular $n$-by-$n$ matrices over a ring $R$. Throughout this note, we will assume the reader has some familiarity with Thompson's groups; see \cite{cannon96} for a standard reference. (As a remark, we do not discuss Thompson's group $T$ in this framework, as it turns out it is somewhat different than $F$ and $V$ from this point of view.)

The Thompson-like group $\Thomp{G_*}$ can be viewed as a sort of limit of the family $(G_n)$. One can compare and contrast it to other limiting operations, e.g., the direct limit. With a direct limit, finiteness properties tend to be destroyed. For example, the symmetric groups $S_n$ are all finitely presented, but their direct limit $S_\infty$ is not even finitely generated. However, arranging the $S_n$ in a natural cloning system and taking the ``Thompson limit'', one gets Thompson's group $V$, which is still finitely presented. This is another motivating factor in the axiomatization of cloning systems; that they yield a limiting procedure that tends to preserve finiteness properties.

Finally, cloning systems, and the Thompson-like groups they produce, simply serve as new examples of interesting groups. For example, in an REU run by Dan Farley, it was shown that certain such examples are \emph{coCF groups} \cite{berns-zieve14}. There is an open conjecture that every coCF group embeds into Thompson's group $V$, and these are potential counterexamples. It remains open whether or not these groups can in fact embed into $V$, but at least there seems to be no natural embedding.

In this note we will first discuss the definition of cloning system in Section~\ref{sec:def}, and discuss various examples in Sections~\ref{sec:ex} and~\ref{sec:new_ex}, then explain how Thompson-like groups and natural cube complexes arise from cloning systems in Sections~\ref{sec:groups} and~\ref{sec:cpx}, and finally discuss the finiteness properties of the Thompson-like groups in Section~\ref{sec:fin}. We will not give any proofs in this expository note (except for proving some statements about the new examples in Section~\ref{sec:new_ex}), but the interested reader can reference \cite{witzel17} for more details.

\subsection*{Acknowledgments} Thanks to Stefan Witzel for suggesting many good improvements to this note, and to an anonymous referee for valuable feedback.


\section{Definitions}\label{sec:def}

This section is devoted to defining a \emph{cloning system} on a family of groups. We fix a family of groups $(G_n)_{n\in\N}$. The rest of the data consist of three families of maps.

\subsection{Directed system morphisms} First, we want the $G_n$ to form a directed system of groups. That is, there should exist maps $\iota_{m,n} \colon G_m \to G_n$, for each $m\le n$, such that $\iota_{n,n}=\id_{G_n}$ for all $n$ and $\iota_{\ell,m}\circ\iota_{m,n}=\iota_{\ell,n}$ for all $\ell\le m\le n$. The astute reader will notice that we are writing composition as though our maps take inputs on the left; indeed we write the $\iota_{m,n}$ maps on the right of their arguments, so notation like $(g)\iota_{m,n}$ is our convention. (This is not really important for this note, but we maintain this convention to be consistent with \cite{witzel17}.) We also require that the $\iota_{m,n}$ be injective, so we can view the direct limit $\colim G_n$ as being the direct union of its subgroups $G_n$. An easy example is the family of symmetric groups $S_n$, with $\iota_{m,n} \colon S_m \to S_n$ given by inclusion, i.e., the image is the subgroup fixing $\{m+1,\dots,n\}$ pointwise.

\subsection{Representation maps} Next, we want to fix a homomorphism $\rho_n \colon G_n \to S_n$ for each $n$. This should be viewed as specifying a way that each $G_n$ acts on the set $\{1,\dots,n\}$, and for the sake of giving $\rho_n$ a name we will call it a \emph{representation map}. For example if $G_n$ is the braid group $B_n$ on $n$ strands, then elements of $B_n$ naturally permute the numbering of the strands, yielding the desired map $B_n\to S_n$. For some choices of $G_n$ there will not be any particularly interesting maps to $S_n$, and in practice the $\rho_n$ will often just be the trivial maps. This is called the \emph{pure} case, and it still yields interesting cloning systems, so a lack of maps to $S_n$ is not a roadblock to finding a cloning system on a given family of groups. We write the maps $\rho_n$ on the left, so notation like $\rho_n(g)$ is our convention. For the $\rho_n$ to count as representation maps, we impose a restriction, namely that the $\rho_n$ should give a \emph{homomorphism of directed systems} $\rho_* \colon G_* \to S_*$. This just means we require
\begin{eqnarray}\label{eq:iota_rho}
 \rho_n((g)\iota_{m,n}) = (\rho_m(g))\iota_{m,n}
\end{eqnarray}
for all $m\le n$ and all $g\in G_m$. (Here we are writing $\iota_{m,n}$ for both the map $G_m\to G_n$ and also the map $S_m\to S_n$.)

\subsection{Cloning maps} Finally, we need our ``cloning maps''. More precisely, we want a family of injective maps from $G_n$ to $G_{n+1}$, for each $n$. There should be $n$ such maps, denoted $\clone_k^n \colon G_n \to G_{n+1}$, for $1\le k\le n$. We mention now and will reiterate later that the $\clone_k^n$ \emph{need not be group homomorphisms}. They are merely (injective) functions on sets. Shortly, when we state the cloning axioms, there will obviously be restrictions on what the $\clone_k^n$ can be, but for now we just have $n$ functions from $G_n$ to $G_{n+1}$. We do impose one restriction: to be termed \emph{cloning maps}, they should satisfy the rule that
\begin{eqnarray}\label{eq:iota_clone}
 \iota_{m,n} \circ \clone_k^n = \clone_k^m \circ \iota_{m+1,n+1}
\end{eqnarray}
for all $1\le k\le m\le n$. Again, the way these equations are written, it is clear that we must write the functions $\clone_k^n$ on the right of their inputs, so notation like $(g)\clone_k^n$ is our convention.

\subsection{The cloning axioms} We now state the axioms for the quadruple
$$((G_n)_{n\in\N},(\iota_{m,n})_{m\le n},(\rho_n)_{n\in\N},(\clone_k^n)_{1\le k\le n})$$
to be a \emph{cloning system}. In the axioms, we always have $1\le k<\ell\le n$ and $g,h\in G_n$.

\textbf{(C1):} (\emph{Cloning a product}) $(gh)\clone_k^n = (g)\clone_{\rho_n(h)k}^n (h)\clone_k^n$.

\textbf{(C2):} (\emph{Product of clonings}) $\clone_\ell^n \circ \clone_k^{n+1} = \clone_k^n \circ \clone_{\ell+1}^{n+1}$.

\textbf{(C3):} (\emph{Compatibility}) $\rho_{n+1}((g)\clone_k^n)(i) = (\rho_n(g))\symmclone_k^n(i)$ for all $i\ne k,k+1$.

\begin{definition}[Cloning system]\label{def:clone}
 Let $((G_n)_{n\in\N},(\iota_{m,n})_{m\le n})$ be a directed system of groups, let $(\rho_n)_{n\in\N}$ be a family of representation maps on the directed system (so Equation~\eqref{eq:iota_rho} is satisfied), and let $((\clone_k^n)_{1\le k\le n})$ be a family of cloning maps on the directed system (so Equation~\eqref{eq:iota_clone} is satisfied). If the quadruple $((G_n)_{n\in\N},(\iota_{m,n})_{m\le n},(\rho_n)_{n\in\N},(\clone_k^n)_{1\le k\le n})$ satisfies (C1), (C2) and (C3) then we call it a \emph{cloning system}.
\end{definition}

In axiom (C3) there is the mysterious notation $\symmclone_k^n$. This will be explained in the section on examples below, see Example~\ref{ex:symm}. (Very quickly: $\symmclone_k^n$ are the cloning maps for the natural cloning system on the symmetric groups.)

We mention some heuristic ways of understanding the axioms. Note that axiom (C1) is saying that the $\clone_k^n$ are not necessarily homomorphisms, but are sort of ``twisted'' homomorphisms, with the twisting given by $\rho_n$. Axiom (C2) is sort of a statement about cloning maps commuting, though the subscripts change in a certain natural way, reminiscent of standard relations in Thompson's group $F$. Finally axiom (C3) says that when hitting everything with the $\rho_n$, the cloning system resembles the standard cloning system on the symmetric groups (which, again, we have not stated yet, but will do so in Example~\ref{ex:symm}). As a remark, in practice (C3) often holds even for $i=k,k+1$, for example in the standard cloning system on the symmetric groups, but this is not axiomatically required. We will discuss some new examples in Section~\ref{sec:new_ex} where (C3) does not hold for $i=k,k+1$.

\medskip

To close out this section on the definition of cloning systems, we discuss one additional property that we usually want, to ensure that a cloning system is ``nice''. This is the property of a cloning system being \emph{properly graded}.

\begin{definition}[Properly graded]\label{def:prop_grad}
 Let $((G_n)_{n\in\N},(\iota_{m,n})_{m\le n},(\rho_n)_{n\in\N},(\clone_k^n)_{1\le k\le n})$ be a cloning system. We call it \emph{properly graded} if for all $1\le k\le n$ we have the inclusion
 $$\image \clone_k^n \cap \image \iota_{n,n+1} \subseteq \image(\iota_{n-1,n} \circ \clone_k^n) \text{.}$$
\end{definition}

In other words, if an element $g \in G_{n+1}$ can be ``uncloned'', and also happens to already lie in $G_n \le G_{n+1}$, then when $g$ is treated as an element of $G_n$ it can still be ``uncloned''. The terminology comes from viewing the $G_n$ as being a filtration, or grading, of $\colim G_n$, and requiring the cloning maps to be well behaved with respect to this grading. The precise nature of the ``niceness'' that cloning systems enjoy when they are properly graded does not really come into play until one builds a Thompson-like group and a Stein--Farley cube complex on which the group acts, after which the groups $G_n$ will appear as vertex stabilizers. This might fail if the cloning system is not properly graded, namely, the stabilizers might not equal the $G_n$. This will all be discussed later, when we construct the groups and complexes.


\section{Existing examples}\label{sec:ex}

We now give some examples of cloning systems. The examples in this section were all given in \cite{witzel17}. In Section~\ref{sec:groups}, we will discuss Thompson-like groups that arise from cloning systems, and in particular we will refer back to these examples to discuss the groups in these cases.

\begin{example}[Direct powers]\label{ex:powers}
 These examples were also observed by Slobodan Tanusevski, and appear in his PhD thesis \cite{tanusevski14}. In discussions with him, we discovered that his construction was an example of a cloning system, albeit in very different language. We state the example here in our language.

 Let $G$ be any group, and consider the family $(G^n)_{n\in\N}$ of direct powers of $G$. The maps $\iota_{m,n} \colon G^m \to G^n$ are given by sending $(g_1,\dots,g_m)$ to $(g_1,\dots,g_m,1,\dots,1)$, where the identity $1\in G$ fills the last $n-m$ entries. The representation maps $\rho_n \colon G^n \to S_n$ are trivial. The cloning map $\clone_k^n \colon G^n \to G^{n+1}$ is given by
 $$(g_1,\dots,g_n)\clone_k^n \defeq (g_1,\dots,g_k,g_k,\dots,g_n) \text{,}$$
 so it copies the $k$th entry into the $(k+1)$st entry and shifts all the higher entries up one, ending with $g_n$ as the $(n+1)$st entry.

 It is a straightforward exercise to check that the quadruple
 $$((G^n)_{n\in\N},(\iota_{m,n})_{m\le n},(\rho_n)_{n\in\N},(\clone_k^n)_{1\le k\le n})$$
 is a cloning system.

 A variation on this example is the following. Let $\phi_1,\phi_2 \colon G \to G$ be two self-monomorphisms of $G$. Then the maps
 $$\clone_k^n \colon (g_1,\dots,g_n) \mapsto (g_1,\dots,\phi_1(g_k),\phi_2(g_k),\dots,g_n)$$
 are still injective, and it is easy to check that they still yield a cloning system on $(G^n)$. In the case when $G$ is finite, $\phi_1=\id_G$ and $\phi_2$ is an automorphism, the resulting Thompson-like groups are known to be coCF \cite{berns-zieve14}, and serve as potential counterexamples to the conjecture that $V$ is universal coCF; see \cite{berns-zieve14} for more details.
\end{example}

\begin{example}[Symmetric groups]\label{ex:symm}
 Take $G_n=S_n$ and let $\iota_{m,n} \colon S_m \to S_n$ be the usual inclusion, whereby we view $S_m$ as the subgroup of $S_n$ fixing $\{m+1,\dots,n\}$ pointwise. Also, naturally, the maps $\rho_n \colon S_n \to S_n$ are just the identity. It remains to state the cloning maps $\symmclone_k^n \colon S_n \to S_{n+1}$ (note the special notation for these cloning maps, which are themselves referenced in axiom (C3)).
 
 Given a permutation $g\in S_n$, the easiest way to understand $(g)\symmclone_k^n \in S_{n+1}$ is with a picture. Since $g$ is a bijection from $\{1,\dots,n\}$ to itself, we can draw $g$ as a diagram of arrows from one copy of $\{1,\dots,n\}$ up to a second copy. Then $(g)\symmclone_k^n$ is the diagram of arrows from $\{1,\dots,n+1\}$ up to a second copy of itself obtained by bifurcating the arrow starting at $k$ into two parallel arrows. See Figure~\ref{fig:symm_clone} for an example.

 \begin{figure}[htb]
 \centering
 \begin{tikzpicture}[line width=0.8pt]
  \draw[->] (0,-2) -- (1,0); \draw[->] (1,-2) -- (0,0);
  \node at (2,-1) {$\stackrel{\symmclone_2^2}{\longrightarrow}$};
  \node at (0,-2.25) {$1$};   \node at (0,0.25) {$1$};   \node at (1,-2.25) {$2$};   \node at (1,0.25) {$2$};
  \begin{scope}[xshift=3cm]
   \draw[->] (0,-2) -- (2,0); \draw (1,-2)[->] -- (0,0); \draw (2,-2)[->] -- (1,0);
   \node at (0,-2.25) {$1$};   \node at (0,0.25) {$1$};   \node at (1,-2.25) {$2$};   \node at (1,0.25) {$2$};   \node at (2,-2.25) {$3$};   \node at (2,0.25) {$3$};
  \end{scope}
 \end{tikzpicture}
 \caption{An example of cloning in symmetric groups. Here we see that $(1~2)\symmclone_2^2 = (1~3~2)$.}
 \label{fig:symm_clone}
 \end{figure}
 
 A more formal definition of $\symmclone_k^n \colon S_n \to S_{n+1}$ is as follows. Let $g\in S_n$, so $g$ is a bijection from $\{1,\dots,n\}$ to itself. We want to specify what $(g)\symmclone_k^n$ is, as a bijection from $\{1,\dots,n+1\}$ to itself. The technical definition is as follows:
 
 $$((g)\symmclone_k^n)m \defeq \left\{\begin{array}{ll}
    gm & \text{if } m\le k \text{ and } gm\le gk \\
    (gm)+1 & \text{if } m<k \text{ and } gm>gk \\
    g(m-1) & \text{if } m>k \text{ and } g(m-1)<gk \\
    g(m-1)+1 & \text{if } m>k \text{ and } g(m-1)\ge gk
    \end{array}\right.$$
 
 For example, in Figure~\ref{fig:symm_clone}, where $g=(1~2)$ and $k=2$, when $m=1$ we have $m\le k$ and $gm=2>1=gk$, so the definition says $((1~2)\symmclone_2^2)1 = ((1~2)1)+1=3$, and indeed the picture of $(1~2)\symmclone_2^2$ shows $1$ going to $3$. As another example, take $m=3$: then $m>k$ and $g(m-1)=1=gk$, so the definition says $((1~2)\symmclone_2^2)3=(1~2)(3-1)+1=2$, and indeed, we see that $3$ goes to $2$.
 
 Using this technical definition, it is possible to formally check that
 $$((G_n)_{n\in\N},(\iota_{m,n})_{m\le n},(\rho_n)_{n\in\N},(\clone_k^n)_{1\le k\le n})$$
 is a properly graded cloning system \cite[Examples~2.9 and~2.16]{witzel17}. As some foreshadowing to Section~\ref{sec:groups}, the Thompson-like group arising from this cloning system will be Thompson's group $V$.
\end{example}

\begin{example}[Braid groups]\label{ex:braids}
 The other example of a cloning system that predates our axiomatization comes from Brin and Dehornoy's \emph{braided Thompson group} $\Vbr$ (often denoted $BV$) \cite{brin07,dehornoy06}.
 
 Take $(G_n)$ to be the family of braid groups, $G_n=B_n$. The inclusion $\iota_{m,n} \colon B_m \to B_n$ is given by adding $n-m$ extra strands to the right of an $m$-strand braid, to get an $n$-strand braid. The representation map $\rho_n \colon B_n \to S_n$ is the usual map taking the numbering of strands at the bottom to the numbering at the top. Finally, the cloning map $\clone_k^n \colon B_n \to B_{n+1}$ is obtained by bifurcating the $k$th strand (counting at the bottom) into two parallel strands, such that no other strands pass between them; see Figure~\ref{fig:braid_clone} for an example.
 
 \begin{figure}[htb]
 \centering
 \begin{tikzpicture}[line width=0.8pt]
  \draw
   (3,-0.5) to [out=-90, in=90, looseness=1] (4,-3);
  \draw[white, line width=4pt]
   (4,-0.5) to [out=-90, in=90, looseness=1] (3,-3);
  \draw
   (4,-0.5) to [out=-90, in=90, looseness=1] (3,-3);
  \node at (5,-1.5) {$\stackrel{\clone_1^2}{\longrightarrow}$};
  \begin{scope}[xshift=3cm]
  \draw
   (3,-0.5) to [out=-90, in=90, looseness=1] (5,-3);
  \draw[white, line width=4pt]
   (5,-0.5) to [out=-90, in=90, looseness=1] (4,-3);
  \draw
   (5,-0.5) to [out=-90, in=90, looseness=1] (4,-3);
  \draw[white, line width=4pt]
   (4,-0.5) to [out=-90, in=90, looseness=1] (3,-3);
  \draw
   (4,-0.5) to [out=-90, in=90, looseness=1] (3,-3);
  \end{scope}
 \end{tikzpicture}\caption{An example of cloning in braid groups. Here we start with a braid $b\in B_2$, apply the cloning map $\clone_1^2$, and get a braid in $B_3$ that looks like $b$ with its first strand cloned (counting at the bottom).}\label{fig:braid_clone}
 \end{figure}
 
 The work involved in checking that $((B_n)_{n\in\N},(\iota_{m,n})_{m\le n},(\rho_n)_{n\in\N},(\clone_k^n)_{1\le k\le n})$ is a properly graded cloning system is similar to that in the symmetric group example. Also, one can instead use the family $(PB_n)$ of pure braid groups, and get a similar cloning system (in this case the representation maps will even be trivial).
\end{example}

\begin{example}[Upper triangular matrix groups]\label{ex:borel}
 In \cite{witzel17}, one of the main new examples of cloning systems involved upper triangular matrix groups. For $R$ a unital ring, let $B_n(R)$ be the group of invertible upper triangular $n$-by-$n$ matrices. Let us describe a properly graded cloning system on the family of groups $(B_n(R))_{n\in\N}$. First, the map $\iota_{m,n} \colon B_m(R) \to B_n(R)$ is the usual one, giving by sending an $m$-by-$m$ matrix $A$ to the $n$-by-$n$ matrix
 $$\begin{pmatrix}
    \begin{array}{cc}
     A & 0 \\ 0 & I_{n-m}
    \end{array}
   \end{pmatrix}\text{.}
 $$
 Next, the representation maps $\rho_n \colon B_n(R) \to S_n$ are all taken to be the trivial map. Finally, the cloning map $\clone_k^n \colon B_n(R) \to B_{n+1}(R)$ is given by
 \[
 \left(
 \begin{array}{ccc}
 A_{<,<} & A_{<,k} & A_{<,>}\\
 0 & A_{k,k} & A_{k,>}\\
 0 & 0 & A_{>,>}
 \end{array}
 \right)\clone_k^n
 =
 \left(
 \begin{array}{cccc}
 A_{<,<} & A_{<,k} & A_{<,k} & A_{<,>}\\
 0 & A_{k,k} & 0 & 0\\
 0 & 0 & A_{k,k} & A_{k,>}\\
 0 & 0 & 0  & A_{>,>}
 \end{array}
 \right)\text{.}
 \]
 
 Here the $A_{*,*}$ represent blocks whose entries lie in positions relative to $k$ as indicated by the subscripts. Note that the block $A_{<,k}$ has width $1$ and the block $A_{k,>}$ has height $1$. It is easier to see what $\clone_k^n$ does by looking at an example. In this example we see what $\clone_3^5$ does to an illustrative $5$-by-$5$ matrix.
 $$\begin{pmatrix} 1&2&3&4&5\\0&6&7&8&9\\0&0&10&11&12\\0&0&0&13&14\\0&0&0&0&15
    \end{pmatrix}\clone_3^5 = \begin{pmatrix}
    1&2&3&3&4&5\\0&6&7&7&8&9\\0&0&10&0&0&0\\0&0&0&10&11&12\\0&0&0&0&13&14\\0&0&0&0&0&15
    \end{pmatrix}$$
 Showing that these data define a properly graded cloning system is not too difficult. The only step that requires some work is axiom (C1), which, since the $\rho_n$ are trivial, says that the cloning maps should be homomorphisms, but technically this ``just'' requires doing a matrix multiplication.
\end{example}


\section{New Examples}\label{sec:new_ex}

Here we present some new examples, which did not appear in \cite{witzel17} but are natural additions to the list of established cloning systems.

\begin{example}[Signed symmetric groups]\label{ex:ssgs}
 Let $G_n=S_n^\pm$ be the \emph{signed symmetric groups}. The group $S_n^\pm$ is the Coxeter group of type $B_n=C_n$, with presentation
 $$S_n^\pm \defeq \left\langle s_1,\dots,s_n \left| \begin{array}{cc}
                                                  s_i^2 = 1 & \text{ for all } i\\
                                                  s_i s_j = s_j s_i & \text{ for } |i-j|>1\\
                                                  (s_i s_{i+1})^3 = 1 & \text{ for } 1\le i\le n-2\\
                                                  (s_{n-1} s_n)^4 = 1
                                                 \end{array}\right.\right\rangle\text{.}$$
 We can realize $S_n^\pm$ as the group of permutations $\sigma$ of $\{1,-1,2,-2,\dots,n,-n\}$ satisfying $\sigma(-i)=-\sigma(i)$ for all $i$. The generators $s_i$ for $1\le i \le n-1$ are the permutations $(i~i+1)(-i~-(i+1))$, and these generate the copy of $S_n$ in $S_n^\pm$ consisting of those $\sigma$ stabilizing the subset $\{1,\dots,n\}$. The generator $s_n$ is the transposition $(n~(-n))$.
 
 We will now put a cloning system on the family $(S_n^\pm)$. We have inclusions $\iota_{m,n} \colon S_m^\pm \to S_n^\pm$ for $m<n$, given by viewing $S_m^\pm$ as the subgroup fixing $\{m+1,-(m+1),\dots,n,-n\}$ pointwise. We also have natural representation maps $\rho_n \colon S_n^\pm \to S_n$ given by sending each $s_i$ for $i<n$ to the generator also called $s_i$ in $S_n$, and sending $s_n$ to the identity. This is clearly a well defined epimorphism, and in fact yields a splitting $S_n^\pm \cong (\Z/2\Z) \wr S_n$, where recall that the wreath product is $(\Z/2\Z) \wr S_n = (\Z/2\Z)^n \rtimes S_n$ with the natural action of $S_n$ on $(\Z/2\Z)^n$. We will not really use this splitting, but it is good to have in mind. We now construct cloning maps $\clone_k^n$ and explain why all the cloning axioms hold.
 
 Since we have an explicit presentation, we will first define the $\clone_k^n \colon S_n^\pm \to S_{n+1}^\pm$ on generators. We declare:
 $$(s_i)\clone_k^n \defeq \left\{\begin{array}{ll}
                                  s_{i+1} & \text{if } k<i<n\\
                                  s_i s_{i+1} & \text{if } k=i<n\\
                                  s_{i+1}s_i & \text{if } k=i+1\le n\\
                                  s_i & \text{if } i+1<k\le n\\
                                  s_{n+1} & \text{if } k<i=n\\
                                  s_{n+1} s_n s_{n+1} & \text{if } k=i=n \text{.}
                                 \end{array}\right.$$
 In the $i<n$ cases these are all the same as for the standard cloning maps on $S_n$. For the last two cases, intuitively, if we view $S_n^\pm$ via pictures like in Figure~\ref{fig:symm_clone}, but now the $i$th arrow has a positive or negative orientation to indicate whether it takes $(i,-i)$ to $(j,-j)$ or $(-j,j)$ for whichever $j$ is appropriate, then $s_n$ looks like the identity except the $n$th arrow gets twisted to the opposite orientation. This makes the second to last case clear, and for the last case note that when we bifurcate this last arrow, it becomes two arrows that cross and both switch orientation, which corresponds to $s_n s_{n-1} s_n$.
 
 To prove this yields a cloning system, we follow the procedure from \cite[Example~9.1]{witzel17}. We need to verify that (C2) and (C3) hold on the generators $s_i$, and then to extend $\clone_k^n$ to be defined on all of $S_n^\pm$ we need to check that if we use (C1) to define $\clone_k^n$ on products of generators, this is well defined according to the defining relations above. The first thing, that (C2) holds on generators, is easy but tedious, and amounts to checking lots of cases, so we leave this to the reader. Checking (C3) on generators is more interesting. For $s_1$ through $s_{n-1}$ it is evident, since $\rho_n$ restricts to the identity on the subgroup $S_n \le S_n^\pm$. Now consider $s_n$. We need to show that $\rho_{n+1}((s_n)\clone_k^n)(i) = (\rho_n(s_n))\symmclone_k^n(i)$ for all $i\ne k,k+1$. The right hand side is $i$, since $\rho_n(s_n)=1$. If $k<n$ then the left hand side is $\rho_{n+1}(s_{n+1})(i)=\id(i)=i$, so this case is done. If $k=n$ then the left hand side is $\rho_{n+1}((s_n)\clone_n^n)(i) = \rho_{n+1}(s_{n+1} s_n s_{n+1})(i) = s_n(i)$, and although $s_n$ is not the identity in $S_{n+1}$, it does send $i$ to $i$ for $i<n$, which is sufficient for (C3) to hold. As a remark, in all the examples in \cite{witzel17}, (C3) held for all $i$, but as this example illustrates we can also have cloning systems where (C3) ``strictly'' holds.
 
 The last thing to check is that the relations are respected upon extending $\clone_k^n$ to all of $S_n^\pm$. The only relation that does not work as in the standard $S_n$ cloning system is the last one, that $(s_{n-1} s_n)^4 = 1$. When we use (C1) to define $\clone_k^n$ on $(s_{n-1} s_n)^4$, we get
 $$((s_{n-1} s_n)^4)\clone_k^n = (s)\clone_{s(k)}^n (t)\clone_{s(k)}^n (s)\clone_k^n (t)\clone_k^n (s)\clone_{s(k)}^n (t)\clone_{s(k)}^n (s)\clone_k^n (t)\clone_k^n$$
 where we write $s\defeq s_{n-1}$ and $t\defeq s_n$ in $S_n^\pm$ for brevity. If $k<n-1$ then this becomes $(s_n s_{n+1})^4$ in $S_{n+1}^\pm$, which is the identity as desired. Now suppose $k=n-1$. Then since $s(n-1)=n$ this becomes $(s)\clone_n^n (t)\clone_n^n (s)\clone_{n-1}^n (t)\clone_{n-1}^n (s)\clone_n^n (t)\clone_n^n (s)\clone_{n-1}^n (t)\clone_{n-1}^n$, which is
 $$s_n s_{n-1} s_{n+1} s_n s_{n+1} s_{n-1} s_n s_{n+1} s_n s_{n-1} s_{n+1} s_n s_{n+1} s_{n-1} s_n s_{n+1}$$
 in $S_{n+1}^\pm$, and applying the relations from $S_{n+1}^\pm$ this becomes
 \begin{align*}
  &s_n s_{n-1} s_{n+1} s_n s_{n-1} s_{n+1} s_n s_{n+1} s_n s_{n+1} s_{n-1} s_n s_{n+1} s_{n-1} s_n s_{n+1} \\
  &= s_n s_{n-1} s_{n+1} s_n s_{n-1} s_n s_{n+1} s_n s_{n-1} s_n s_{n+1} s_{n-1} s_n s_{n+1} \\
  &= s_n s_{n+1} s_n s_{n-1} s_{n+1} s_n s_{n-1} s_n s_{n+1} s_{n-1} s_n s_{n+1} \\
  &= s_n s_{n+1} s_n s_{n+1} s_n s_{n+1} s_n s_{n+1} = 1 \text{.}
 \end{align*}
 Finally, when $k=n$ we get a similar expression and it similarly becomes $1$ in $S_{n+1}^\pm$.
 
 We conclude that $((S_n^\pm)_{n\in\N},(\iota_{m,n})_{m\le n},(\rho_n)_{n\in\N},(\clone_k^n)_{1\le k\le n})$ is a cloning system. It is also easy to check that it is properly graded, like in the standard $S_n$ case.
\end{example}

\begin{example}[Twisted braid groups]\label{ex:tbgs}
 One can view $S_n^\pm$ as a ``twisted'' symmetric group, with the generator called $s_n$ in the previous example serving to twist the $n$th arrow, and conjugates of $s_n$ twisting the other arrows. (In this case the twisting move has order $2$, so ``flipping'' might be a better word.) Passing from symmetric groups to braid groups, it makes sense to again try and ``twist'' things. A first attempt at twisted braid groups could be $B_n^\pm = (\Z/2\Z) \wr B_n$, analogous to how $S_n^\pm = (\Z/2\Z) \wr S_n$, with the action of $B_n$ on $(\Z/2\Z)^n$ given via the projection $B_n \to S_n$. This admits a cloning system in much the same way as the previous example, the details of which we will leave to the reader.

More interesting is to view the strands in the braids as ribbons that can twist, and consider the group $B_n^{twist} \defeq \Z \wr B_n = \Z^n \rtimes B_n$, with the action of $B_n$ on $\Z^n$ given via the projection $B_n \to S_n$. Now a twist is an infinite-order move, and it is straightforward to see that ribbon pictures correctly model this wreath product. Since $B_n^{twist}$ is a wreath product, we can write down a presentation. Write the generators as $s_1,\dots,s_{n-1},s_n$, where $s_1,\dots,s_{n-1}$ are the standard generators of $B_n$, and $s_n$ is a twist of the $n$th ribbon. The defining relations are the usual braid relations $s_i s_{i+1} s_i = s_{i+1} s_i s_{i+1}$ for $1\le i\le n-2$ plus the relations $s_n^{-1} s_i^{-1} s_n^{-1} s_i s_n s_i^{-1} s_n s_i = 1$ for all $i$ (to ensure $s_n$ commutes with all its conjugates).

We can define a cloning system on $B_n^{twist}$ as follows. First define $\clone_k^n$ on the generators exactly as in Example~\ref{ex:ssgs}, via
$$(s_i)\clone_k^n \defeq \left\{\begin{array}{ll}
                                  s_{i+1} & \text{if } k<i<n\\
                                  s_i s_{i+1} & \text{if } k=i<n\\
                                  s_{i+1}s_i & \text{if } k=i+1\le n\\
                                  s_i & \text{if } i+1<k\le n\\
                                  s_{n+1} & \text{if } k<i=n\\
                                  s_{n+1} s_n s_{n+1} & \text{if } k=i=n \text{.}
                                 \end{array}\right.$$
																
See Figure~\ref{fig:twist_clone} for an example of this last cloning move. The figure depicts the cloning move $(s_1)\clone_1^1 = s_2 s_1 s_2$ (drawn here as $s_1(s_1^{-1}s_2s_1)s_2$).

\begin{figure}[htb]
 \centering
 \begin{tikzpicture}[line width=0.8pt]
  \draw
	 (0,0) -- (1,0)
   (1,0) to [out=-90, in=90, looseness=1] (0,-2);
	\draw[white,line width=4pt]
	 (0,0) to [out=-90, in=90, looseness=1] (1,-2);
	\draw
	 (0,0) to [out=-90, in=90, looseness=1] (1,-2)
	 (0,-2) to [out=-90, in=90, looseness=1] (-1,-4)
	 (1,-2) to [out=-90, in=90, looseness=1] (2,-4)
	 (0.5,-3) to [out=-180, in=90, looseness=1] (0,-4)
	 (0.5,-3) to [out=0, in=90, looseness=1] (1,-4)
	 (-1,-4) -- (0,-4)   (1,-4) -- (2,-4);
  \node at (3,-1.5) {$\stackrel{\clone_1^1}{\longrightarrow}$};
  \begin{scope}[xshift=5cm]
   \draw
    (0,0) -- (1,0)
		(0,0) to [out=-90, in=90, looseness=1] (-1,-1.5)
	  (1,0) to [out=-90, in=90, looseness=1] (2,-1.5)
	  (0.5,-1) to [out=-180, in=90, looseness=1] (0,-1.5)
	  (0.5,-1) to [out=0, in=90, looseness=1] (1,-1.5);
	 \draw
	  (1,-1.5) to [out=-90, in=90, looseness=1] (-1,-3.5)
		(2,-1.5) to [out=-90, in=90, looseness=1] (0,-3.5);
	 \draw[white,line width=4pt]
	  (-1,-1.5) to [out=-90, in=90, looseness=1] (1,-3.5)
	  (0,-1.5) to [out=-90, in=90, looseness=1] (2,-3.5);
	 \draw[white,line width=10pt]
	  (-.5,-1.5) to [out=-90, in=90, looseness=1] (1.5,-3.5);
	 \draw
	  (-1,-1.5) to [out=-90, in=90, looseness=1] (1,-3.5)
	  (0,-1.5) to [out=-90, in=90, looseness=1] (2,-3.5)
	  (0,-3.5) to [out=-90, in=90, looseness=1] (-1,-4.5)
	  (2,-3.5) to [out=-90, in=90, looseness=1] (1,-4.5);
	 \draw[white,line width=4pt]
	  (-1,-3.5) to [out=-90, in=90, looseness=1] (0,-4.5)
		(1,-3.5) to [out=-90, in=90, looseness=1] (2,-4.5);
	 \draw
	  (-1,-3.5) to [out=-90, in=90, looseness=1] (0,-4.5)
		(1,-3.5) to [out=-90, in=90, looseness=1] (2,-4.5)
		(-1,-4.5) -- (0,-4.5)   (1,-4.5) -- (2,-4.5);
  \end{scope}
 \end{tikzpicture}
 \caption{Cloning a ribbon twist. We see that $(s_1)\clone_1^1 = s_2 s_1 s_2$.}
 \label{fig:twist_clone}
 \end{figure}

For the maps $\rho_n \colon B_n^{twist} \to S_n$ we just take the composition of $B_n^{twist} \to B_n$ with $B_n \to S_n$. Now we can verify the cloning axioms, again following the procedure from \cite[Example~9.1]{witzel17}. Most of the things to check work in exactly the same way as they did in Example~\ref{ex:ssgs}. In fact the only step that does not follow identically is the verification that the relations $s_n^{-1} s_i^{-1} s_n^{-1} s_i s_n s_i^{-1} s_n s_i = 1$ are respected by the cloning maps. We rewrite this as $s_n s_i^{-1} s_n s_i = s_i^{-1} s_n s_i s_n$ and apply $\clone_k^n$ to both sides. On the left we get
$$(s_n s_i^{-1} s_n s_i)\clone_k^n = (s_n)\clone_k^n (s_i^{-1})\clone_{(i~i+1)k}^n (s_n)\clone_{(i~i+1)k}^n (s_i)\clone_k^n$$
and on the right we get
$$(s_i^{-1} s_n s_i s_n)\clone_k^n = (s_i^{-1})\clone_{(i~i+1)k}^n (s_n)\clone_{(i~i+1)k}^n (s_i)\clone_k^n (s_n)\clone_k^n\text{.}$$
It is now an exercise to verify that these equal each other in $B_{n+1}^{twist}$, for each $1\le i\le n-1$ and $1\le k\le n$. We will work out the most difficult case, when $i=n-1$ and $k=n$. In this case we need to show that $(s_n)\clone_n^n (s_{n-1}^{-1})\clone_{n-1}^n (s_n)\clone_{n-1}^n (s_{n-1})\clone_n^n$ equals $(s_{n-1}^{-1})\clone_{n-1}^n (s_n)\clone_{n-1}^n (s_{n-1})\clone_n^n (s_n)\clone_n^n$. Applying all the cloning maps, the thing to show is that $s_{n+1}s_n s_{n+1} (s_n s_{n-1})^{-1} s_{n+1} s_n s_{n-1}$ equals $(s_n s_{n-1})^{-1} s_{n+1} s_n s_{n-1} s_{n+1}s_n s_{n+1}$, or equivalently that $(s_n s_{n-1})^{-1} s_{n+1} s_n s_{n-1}$ commutes with $s_{n+1}s_n s_{n+1}$. Note that the former is in the $\Z^{n+1}$ factor of $B_{n+1}^{twist}=\Z^{n+1}\rtimes B_n$, so conjugating it by $s_{n+1}s_n s_{n+1}$ is the same as conjugating it by $s_n$, which yields
$$(s_n s_{n-1}s_n)^{-1} s_{n+1} s_n s_{n-1}s_n = (s_{n-1} s_n s_{n-1})^{-1} s_{n+1} s_{n-1}s_n s_{n-1} = (s_n s_{n-1})^{-1} s_{n+1} s_n s_{n-1}\text{.}$$
Hence $(s_n s_{n-1})^{-1} s_{n+1} s_n s_{n-1}$ commutes with $s_{n+1}s_n s_{n+1}$ and we are done.

As a remark, one could also construct various ``pure'' versions, restricting to only pure braids and/or ``pure twists''. This is a straightforward generalization that we leave to the reader.
\end{example}


\section{Thompson-like groups}\label{sec:groups}

Given a cloning system $((G_n)_{n\in\N},(\iota_{m,n})_{m\le n},(\rho_n)_{n\in\N},(\clone_k^n)_{1\le k\le n})$, as we have repeatedly promised, one gets a Thompson-like group $\Thomp{G_*}$ out of it. This group contains Thompson's group $F$ and all the groups $G_n$ as subgroups, and it maps to Thompson's group $V$, all in natural ways. In this section we will describe the group, and state some of the important properties. We will not get into the formal details that are necessary to define the group, since this would take us on a long, technical detour into \emph{Brin--Zappa-Sz\'ep products}. These details can be found in Sections~1 and~2 of \cite{witzel17}. Rather than get into these details, here we will simply state what elements of the group $\Thomp{G_*}$ look like, discuss some basic properties of the groups, and give examples.

\subsection{Elements} Throughout, we fix a cloning system $((G_n)_{n\in\N},(\iota_{m,n})_{m\le n},(\rho_n)_{n\in\N},(\clone_k^n)_{1\le k\le n})$. An element of the Thompson-like group $\Thomp{G_*}$ is represented by a triple $(T_-,g,T_+)$. Here $T_-$ and $T_+$ are trees (by which we always mean finite rooted binary trees) with the same number of leaves, say $n$, and $g\in G_n$. This triple \emph{represents} an element of $\Thomp{G_*}$ in the sense that elements of $\Thomp{G_*}$ are actually equivalence classes of such triples, under a certain equivalence relation. The equivalence relation is the symmetric transitive hull of moves called \emph{expansions}. An expansion of $(T_-,g,T_+)$ is a triple $(U_-,h,U_+)$, where $U_+$ is the tree obtained from $T_+$ by adding a caret to the $k$th leaf of $T_+$, for some $1\le k\le n$, $U_-$ is the tree obtained by adding a caret to the $(\rho_n(g)k)$th leaf of $T_-$, and $h=(g)\clone_k^n$. So, two triples are considered \emph{equivalent} if one can get from one to the other by a finite sequence of expansions and inverse expansions (called \emph{reductions}). It is worth pointing out that the representation maps $\rho_n$ and the cloning maps $\clone_k^n$ both came in to play in defining this equivalence relation. See Figure~\ref{fig:borel_thomp_expansion} for an example of expansion in the matrix group case.

\begin{figure}[htb]
 \centering
\begin{tikzpicture}[line width=0.8pt, scale=0.4]
  \begin{scope}
  \node at (-3,-1){$\Bigg($};
  \draw
   (-2,-2) -- (0,0) -- (2,-2)   (0,-2) -- (-1,-1);
  \filldraw
   (-2,-2) circle (1.5pt)   (0,0) circle (1.5pt)   (2,-2) circle (1.5pt)   (-1,-1) circle (1.5pt)   (0,-2) circle (1.5pt);
   \node at (2.5,-2){,};
   \node at (5,-1){$\begin{pmatrix}
                    1&2&3\\0&4&5\\0&0&6
                   \end{pmatrix}$};
                   
   \node at (7.5,-2){,};
   \begin{scope}[xshift=4in]    
    \draw
   (-2,-2) -- (0,0) -- (2,-2)   (0,-2) -- (1,-1);
  \filldraw
   (-2,-2) circle (1.5pt)   (0,0) circle (1.5pt)   (2,-2) circle (1.5pt)   (1,-1) circle (1.5pt)   (0,-2) circle (1.5pt);
   \node at (3,-1){$\Bigg)$};
   \end{scope}
   \end{scope}

    \node at (14.5,-1){$\longrightarrow$};
   
   \begin{scope}[xshift=18.75cm]
    \node at (-3,-1){$\Bigg($};
    \draw
   (-2,-2) -- (0,0) -- (2,-2)   (0,-2) -- (-1,-1)   (-1,-2) -- (-.5,-1.5);
    \filldraw
   (-2,-2) circle (1.5pt)   (0,0) circle (1.5pt)   (2,-2) circle (1.5pt)   (-1,-1) circle (1.5pt)   (0,-2) circle (1.5pt)   (-1,-2) circle (1.5pt);
   \node at (3,-2){,};
   \node at (6.5,-1){$\begin{pmatrix}
                    1&2&2&3\\0&4&0&0\\0&0&4&5\\0&0&0&6
                   \end{pmatrix}$};
                   
   \node at (9.75,-2){,};
   \begin{scope}[xshift=4.8in]    
    \draw
   (-2,-2) -- (0,0) -- (2,-2)   (0,-2) -- (1,-1)   (.5,-1.5) -- (1,-2);
  \filldraw
   (-2,-2) circle (1.5pt)   (0,0) circle (1.5pt)   (2,-2) circle (1.5pt)   (1,-1) circle (1.5pt)   (0,-2) circle (1.5pt)   (1,-2) circle (1.5pt);
   \node at (2.8,-1){$\Bigg)$};
   \end{scope}
   \end{scope}
  
\end{tikzpicture}
\caption[Expansion]{An example of expansion in $\Thomp{B_*(\Q)}$. The expansion amounts to adding a caret to the second leaf of each tree (note that $\rho_3$ is trivial) and applying $\clone_2^3$ to the matrix.}
\label{fig:borel_thomp_expansion}
\end{figure}

We write elements of $\Thomp{G_*}$ as $[T_-,g,T_+]$, to mean the equivalence class of the triple $(T_-,g,T_+)$. We are claiming that $\Thomp{G_*}$ is a group, so it had better make sense to take a product $[T_-,g,T_+][U_-,h,U_+]$. It turns out that the following is the right thing to do. Find a common expansion $S$ of $T_+$ and $U_-$, for instance their union, and then perform expansions on the two triples to get $[T_-,g,T_+]=[T_-',g',S]$ and $[U_-,h,U_+]=[S,h',U_+']$ for some $T_-'$, $g'$, $h'$ and $U_+'$. then the product is defined by
$$[T_-,g,T_+][U_-,h,U_+] \defeq [T_-',g'h',U_+'] \text{.}$$
This turns out to be a well defined group operation. As a remark, expanding $T_+$ to $S$ amounts to adding carets. Hence, the corresponding expansion of $(T_-,g,T_+)$ to $(T_-',g',S)$ amounts to also adding carets to $T_-$ to get $T_-'$ and applying some sequence of cloning maps to $g$ to get $g'$. Similarly $h'$ is just $h$ fed into some sequence of cloning maps. See Figure~\ref{fig:V_mult} for an example of multiplication in $\Thomp{S_*}=V$.

\begin{figure}[ht]
 \centering
\begin{tikzpicture}[line width=0.8pt, scale=0.4]
  \node at (-1.5,-1){$\Big[$};
  \draw (-1,-1) -- (0,0) -- (1,-1);
  \filldraw (-1,-1) circle (1.5pt)   (0,0) circle (1.5pt)   (1,-1) circle (1.5pt);
  \node at (1.5,-1){,};
  \node at (3,-1){$(1~2)$};
  \node at (4.5,-1){,};
   \begin{scope}[xshift=2.35in]
   \draw (-1,-1) -- (0,0) -- (1,-1);
   \filldraw (-1,-1) circle (1.5pt)   (0,0) circle (1.5pt)   (1,-1) circle (1.5pt);
   \node at (1.5,-1){$\Big]$};
   \end{scope}
  \begin{scope}[xshift=4.3in]
   \node at (-3,-1){$\Bigg[$};
   \draw (-2,-2) -- (0,0) -- (2,-2)   (0,-2) -- (1,-1);
   \filldraw (-2,-2) circle (1.5pt)   (0,0) circle (1.5pt)   (2,-2) circle (1.5pt)   (1,-1) circle (1.5pt)   (0,-2) circle (1.5pt);
   \node at (2.5,-2){,};
   \node at (3.25,-1){$\id$};
   \node at (4,-2){,};
   \begin{scope}[xshift=2.5in]
    \draw (-2,-2) -- (0,0) -- (2,-2)   (0,-2) -- (-1,-1);
    \filldraw (-2,-2) circle (1.5pt)   (0,0) circle (1.5pt)   (2,-2) circle (1.5pt)   (-1,-1) circle (1.5pt)   (0,-2) circle (1.5pt);
    \node at (3,-1){$\Bigg]$};
   \end{scope}
  \end{scope}

    \node at (24,-1){$=$};
   
  \begin{scope}[yshift=-2in]
   \node at (-3,-1){$\Bigg[$};
   \draw (-2,-2) -- (0,0) -- (2,-2)   (0,-2) -- (-1,-1);
   \filldraw (-2,-2) circle (1.5pt)   (0,0) circle (1.5pt)   (2,-2) circle (1.5pt)   (-1,-1) circle (1.5pt)   (0,-2) circle (1.5pt);
   \node at (2.5,-2){,};
   \node at (4,-1){$(1~3~2)$};
   \node at (5.75,-2){,};
   \begin{scope}[xshift=3.25in]
    \draw (-2,-2) -- (0,0) -- (2,-2)   (0,-2) -- (1,-1);
    \filldraw (-2,-2) circle (1.5pt)   (0,0) circle (1.5pt)   (2,-2) circle (1.5pt)   (1,-1) circle (1.5pt)   (0,-2) circle (1.5pt);
    \node at (3,-1){$\Bigg]$};
   \end{scope}
   \begin{scope}[xshift=6in]
    \node at (-3,-1){$\Bigg[$};
    \draw (-2,-2) -- (0,0) -- (2,-2)   (0,-2) -- (1,-1);
    \filldraw (-2,-2) circle (1.5pt)   (0,0) circle (1.5pt)   (2,-2) circle (1.5pt)   (1,-1) circle (1.5pt)   (0,-2) circle (1.5pt);
    \node at (2.5,-2){,};
    \node at (3.25,-1){$\id$};
    \node at (4,-2){,};
    \begin{scope}[xshift=2.5in]
     \draw (-2,-2) -- (0,0) -- (2,-2)   (0,-2) -- (-1,-1);
     \filldraw (-2,-2) circle (1.5pt)   (0,0) circle (1.5pt)   (2,-2) circle (1.5pt)   (-1,-1) circle (1.5pt)   (0,-2) circle (1.5pt);
     \node at (3,-1){$\Bigg]$};
    \end{scope}
   \end{scope}
  \end{scope}
  
    \node at (26,-6){$=$};
    
  \begin{scope}[yshift=-4in]
   \node at (-3,-1){$\Bigg[$};
   \draw (-2,-2) -- (0,0) -- (2,-2)   (0,-2) -- (-1,-1);
   \filldraw (-2,-2) circle (1.5pt)   (0,0) circle (1.5pt)   (2,-2) circle (1.5pt)   (-1,-1) circle (1.5pt)   (0,-2) circle (1.5pt);
   \node at (2.5,-2){,};
   \node at (3.75,-1){$(1~3~2)$};
   \node at (5,-2){,};
   \begin{scope}[xshift=3in]
    \draw (-2,-2) -- (0,0) -- (2,-2)   (0,-2) -- (-1,-1);
    \filldraw (-2,-2) circle (1.5pt)   (0,0) circle (1.5pt)   (2,-2) circle (1.5pt)   (-1,-1) circle (1.5pt)   (0,-2) circle (1.5pt);
    \node at (3,-1){$\Bigg]$};
   \end{scope}
  \end{scope}
\end{tikzpicture}
\caption[]{Multiplication in $\Thomp{S_*}=V$. Note that we first have to expand the left triple, using the cloning map $\symmclone_2^2$.}
\label{fig:V_mult}
\end{figure}

\subsection{Basic properties} Section~3 of \cite{witzel17} details some basic properties of the Thompson-like groups $\Thomp{G_*}$, regardless of the cloning system. We list a few of them here. First, the map $[T_-,g,T_+]\mapsto [T_-,\rho_n(g),T_+]$ is a homomorphism $\Thomp{G_*}\to V$ to Thompson's group $V=\Thomp{S_*}$. We call the kernel $\Thkern{G_*}$. If the $\rho_n$ are all the trivial map, i.e., if we are in the pure case, then the image of this map is Thompson's group $F$, and moreover the map splits (otherwise the map does not necessarily split), so in the pure case we have $\Thomp{G_*}=\Thkern{G_*}\rtimes F$ \cite[Observation~3.2]{witzel17}. Second, for $n\in\N$ and $T$ any tree with $n$ leaves, the map $g\mapsto [T,g,T]$ is a monomorphism $G_n \hookrightarrow \Thomp{G_*}$. Hence $\Thomp{G_*}$ always contains the groups $G_n$, for all $n\in\N$, as subgroups \cite[Observation~3.1]{witzel17}. Also, any $\Thomp{G_*}$ also contains Thompson's group $F$ in a natural way, namely as the subgroup of elements of the form $[T_-,1,T_+]$. Finally, if we replace finitely many groups in the family $(G_n)$ with the trivial group, the isomorphism type of $\Thomp{G_*}$ does not change \cite[Proposition~3.6]{witzel17}; this is one way in which we can think of $\Thomp{G_*}$ as a sort of limit of the $G_n$, in that it is immune to changes in an initial segment of the sequence.

\medskip

\subsection{Examples} We quickly review the examples of cloning systems given in Sections~\ref{sec:ex} and~\ref{sec:new_ex}, and discuss some details of the arising Thompson-like groups.

The very first example uses the trivial cloning system, where $G_n=\{1\}$ for all $n$ and all the $\iota_{m,n}$, $\rho_n$ and $\clone_k^n$ are trivial. In this case the Thompson-like group $\Thomp{\{1\}}$ is Thompson's group $F$. Elements are equivalence classes of tree pairs $[T_-,T_+]$.

An element of the Thompson-like group $\Thomp{G^*}$, coming from the family of direct powers $G^n$ of a fixed group $G$ (Example~\ref{ex:powers}), looks like $[T_-,(g_1,\dots,g_n),T_+]$, where $T_-$ and $T_+$ are trees with $n$ leaves. Some facts worth mentioning here are that $\Thomp{G^*}=\Thkern{G^*}\rtimes F$ (since the $\rho_n$ are all trivial), and we actually have a sequence of maps $G\to \Thomp{G^*} \to G$ composing to the identity. Namely, the first map is $g\mapsto [1,g,1]$, where $1$ is the trivial tree, and the second map is $[T_-,(g_1,\dots,g_n),T_+] \mapsto g_1$. These are well defined group homomorphisms, and their composition is clearly the identity on $G$. Hence we have what is called a \emph{retract} of $\Thomp{G^*}$ onto $G$, which has various implications, e.g., for finiteness properties; see \cite{tanusevski14} and \cite[Section~6]{witzel17}.

Next we have Example~\ref{ex:symm} and the cloning system on the family of symmetric groups. Here we find that $\Thomp{S_*}$ equals Thompson's group $V$. Elements are equivalence classes $[T_-,\sigma,T_+]$, where $\sigma$ is a permutation on the numbering of the leaves of the trees. Following this, the cloning system on the (pure) braid groups (Example~\ref{ex:braids}) yields the Thompson-like groups $\Thomp{B_*}=\Vbr$ and $\Thomp{PB_*}=\Fbr$, the braided Thompson groups. The group $\Fbr$ split surjects onto $F$, and we have $\Fbr=\Thkern{PB_*}\rtimes F$. The group $\Vbr$ surjects onto $V$, but this does not split; indeed, $\Vbr$ is torsion-free and $V$ is not. Next, Example~\ref{ex:borel} gives us a cloning system on $(B_n(R))$, and the resulting Thompson-like group $\Thomp{B_*(R)}$ has elements represented as $[T_-,A,T_+]$ for $A$ an invertible upper triangular matrix. Since the representation maps are all trivial in this case, $\Thomp{B_*(R)}$ is a semidirect product of $F$ with $\Thkern{B_*(R)}$.

Moving on to Section~\ref{sec:new_ex}, we first have Example~\ref{ex:ssgs}. The cloning system on $(S_n^\pm)$ gives us a Thompson-like group $V^\pm \defeq \Thomp{S_*^\pm}$. It is not hard to see that $V^\pm$ is an example of the sort of group considered in \cite{bleak17}, namely $V^\pm = V_2(S_2)$ in their notation. Since $S_2$ is semiregular in itself, \cite[Theorem~1]{bleak17} says that $V^\pm \cong V$ (see also the discussion of this precise example in \cite[Section~5]{bleak17}). This is interesting since the families $(S_n^\pm)_n$ and $(S_n)_n$ are different, but yield isomorphic ``Thompson limits''. Finally we have Example~\ref{ex:tbgs} involving the family of twisted braid groups $(B_n^{twist})$, which gives us a Thompson-like group $\Vbr^{twist} \defeq \Thomp{B_*^{twist}}$. The cloning systems on $(B_n^{twist})$ and $(S_n^\pm)$ are compatible with the surjections $B_n^{twist}\to S_n^\pm$, and so we get a surjection $\Vbr^{twist} \to V^\pm \cong V$. It is unclear whether $\Vbr^{twist}$ and $\Vbr$ are isomorphic or not, though we suspect they are not. It is also unclear whether $\Vbr^{twist}$ surjects onto $\Vbr$; the maps $B_n^{twist}\to B_n$ are not compatible with the cloning systems, so there is no ``obvious'' way to build such a surjection.


\section{Stein--Farley complexes}\label{sec:cpx}

The data in a cloning system not only give rise to a Thompson-like group $\Thomp{G_*}$, but also a natural contractible cube complex $\Stein{G_*}$ on which the group acts, called the \emph{Stein--Farley complex}. The action is not cocompact, and is in general not proper, but it nonetheless reveals information about the group. In particular, there is a natural cocompact filtration of the space, and understanding the topology of this filtration reduces via discrete Morse theory to understanding the topology of certain \emph{descending links} of vertices in the space. Also, if the cloning system is properly graded, then the groups $G_n$ are precisely the vertex stabilizers for the action of $\Thomp{G_*}$ on $\Stein{G_*}$. In particular, the situation is ready-made for Brown's Criterion for finiteness properties \cite{brown87}, and one can learn a lot about finiteness properties of $\Thomp{G_*}$ by understanding the aforementioned descending links.

Our goal here is not to get into the details of Morse theory, descending links, or the more complicated situation in \cite[Section~5.5]{witzel17} where one cannot use Brown's Criterion. Rather, in this section we will just state the construction of the Stein--Farley complex and the natural filtration.

\subsection{Vertices} Given a cloning system on a family of groups $(G_n)$, in the previous section we described how to build the Thompson-like group $\Thomp{G_*}$, and there is a similar procedure for building the Stein--Farley complex $\Stein{G_*}$. First, the vertices of the cube complex are equivalence classes $[T,g,E]$ of triples $(T,g,E)$, where $T$ is a tree, say with $n$ leaves, $g$ is an element of $G_n$, and $E$ is a forest with $n$ leaves (and some number of roots). By a \emph{forest} we always mean a finite ordered disjoint union of trees. The next definition will come up a lot and is worth recording now, before stating what equivalence relations gets us from $(T,g,E)$ to $[T,g,E]$.

\begin{definition}[Number of feet]\label{def:feet}
 The \emph{number of feet} of the triple $(T,g,E)$ as above is the number of roots of the forest $E$.
\end{definition}

The equivalence relation we want to impose on such triples $(T,g,E)$ is similar to that for elements of $\Thomp{G_*}$. Namely, we can add appropriate carets to $T$ and $E$ and apply appropriate cloning maps to $g$ without changing the vertex. Details of these moves are exactly the same as for elements of $\Thomp{G_*}$, so we will not repeat them here. There is a new aspect though, which we want to also include in the equivalence relation: $G_n$ acts on the set of triples $(T,g,E)$ with $n$ feet, from the right, and we also mod out this action. Roughly, the action exists due to the fact that if $E$ has $n$ roots then the product $[T,g,E][1_n,h,1_n]$ makes sense, where $1_n$ is the trivial forest consisting of $n$ trivial trees. For each $n$ we mod out the action of $G_n$ on the set of triples with $n$ feet.

\subsection{Edges and cubes} So far we have only defined vertices of $\Stein{G_*}$. We need to say what the cubes of this cube complex are. First, an edge is defined by its endpoints: a ``top'' vertex $[T,g,E]$ and a ``bottom'' vertex $[T,g,E']$, where $E'$ is obtained from $E$ by adding a new caret, whose leaves are the roots of two of the trees in $E$ (so $E$ must have had at least two trees in it). Then, the higher dimensional cubes are glued in anytime the $1$-skeleton of a cube appears. More precisely, a $k$-cube is given by a top vertex and a collection of $k$ \emph{disjoint} new carets that may be added. See Section~4, specifically Subsection~4.3 of \cite{witzel17} for more details, including all the formalism. As seen in Proposition~4.8 of \cite{witzel17}, the space $\Stein{G_*}$ is contractible.

\subsection{The action} Much like one can define a group multiplication operation on two triples $[T_-,g,T_+],[U_-,h,U_+]$, now that we are considering forests and not just trees, we can define a group \emph{action}, namely, the product $[T_-,g,T_+][U,h,E]$ makes sense even if $E$ is a forest (but $U$ should still be a tree). The left triple is a group element and the right triple is a vertex of the space. (As a remark, if $E$ is not a tree, the product $[U,h,E][T_-,g,T_+]$ does not make sense.) The action preserves adjacency of vertices, and so takes cubes to cubes. In summary, the group $\Thomp{G_*}$ acts cellularly on the contractible space $\Stein{G_*}$.

The action of $\Thomp{G_*}$ on $\Stein{G_*}$ is not cocompact, since it leaves invariant the measurement ``number of feet'' on vertices, and this takes infinitely many values, but if we define $\Stein{G_*}^{f\le n}$ to be the full subcomplex supported on those vertices with at most $n$ feet, then it turns out each $\Stein{G_*}^{f\le n}$ is invariant and cocompact \cite[Lemma~5.5]{witzel17}. Hence these subcomplexes provide a cocompact filtration of the cube complex. In the case when the cloning system is properly graded, it turns out that the stabilizer in $\Thomp{G_*}$ of a vertex with $n$ feet is isomorphic to $G_n$ \cite[Lemma~4.9]{witzel17}. Hence if the groups $G_n$ that we started with are ``nice'' in some way, then we now have an action with nice stabilizers.

In particular, having an action on a contractible space with nice stabilizers and a cocompact filtration is the setup of Brown's Criterion for finiteness properties, referenced earlier. Rather than discuss why/whether all this setup tells us things about finiteness properties of $\Thomp{G_*}$, in the last section we will just recall what we mean by finiteness properties, state some results for Thompson-like groups, and discuss the general behavior that seems to often occur.


\section{Finiteness properties}\label{sec:fin}

In this final section, we discuss one of the main subjects that spurred us to develop cloning systems, namely \emph{finiteness properties} of groups. By the finiteness properties of a group, we mean the properties of being of type $\F_n$, for $n\in\N$. We say a group $G$ is a of \emph{type $\F_n$} if it admits a $K(G,1)$ with compact $n$-skeleton. Here a $K(G,1)$, also called a \emph{classifying space} for $G$, is a connected CW complex $X$ with $\pi_1(G)\cong G$ and $\pi_k(X)=0$ for $k\ge2$. For a given $G$, such spaces are unique up to homotopy equivalence. For small $n$, these have nice algebraic interpretations: a group is of type $\F_1$ if and only if it is finitely generated, and type $\F_2$ if and only if it is finitely presented. If a group is of type $\F_n$ for all $n$ we say it is of \emph{type $\F_\infty$}. Every group is of type $\F_0$.

\begin{definition}[Finiteness length]\label{def:fin_length}
 The \emph{finiteness length} $\phi(G)$ of a group $G$ is the largest $n\in\N_0\cup\{\infty\}$ such that $G$ is of type $\F_n$.
\end{definition}

For example, if $G$ is not finitely generated then $\phi(G)=0$, and if $G$ is finitely generated but not finitely presented then $\phi(G)=1$. This finiteness length function $\phi$ is in general not well behaved with respect to limiting procedures. Two standard notions of ``limit'' in group theory are the direct limit (of groups in a directed system) and the inverse limit (of groups in an inverse system). For both of these limiting processes, $\phi$ is poorly behaved. Namely, outside some trivial cases, a direct or inverse limit of infinitely many finite groups will not even be finitely generated, so $\phi(\lim G_n)=0$ even though $\phi(G_n)=\infty$ for all $n$. This can be viewed as a strong failure of $\phi$ to be ``continuous''.

In contrast, the ``Thompson limit'' $\Thomp{G_*}$ of a family of groups $G_n$ with a cloning system seems to be well behaved with respect to $\phi$. The behavior that we find in many examples is:
\begin{eqnarray}\label{eq:fin}
 \phi(\Thomp{G_*}) = \liminf_{n\to\infty}\phi(G_n) \text{.}
\end{eqnarray}
The examples where Equation~\eqref{eq:fin} is known to hold include the following groups. First, when each $G_n=\{1\}$ is trivial, the right hand side of Equation~\eqref{eq:fin} is $\infty$, and $\Thomp{\{1\}}=F$ is of type $\F_\infty$ \cite{brown84} so the left hand side is also $\infty$. Similarly, when $G_n=S_n$, the right hand side is $\infty$ and so is the left since $\Thomp{S_*}=V$ is of type $\F_\infty$ \cite{brown87}. In the braided cases, $B_n$ and $PB_n$ are $\F_\infty$, and so are $\Vbr$ and $\Fbr$ \cite{bux16}, so again Equation~\eqref{eq:fin} reads $\infty=\infty$. When $G_n=G^n$, the $n$th direct power of $G$, it is known \cite{tanusevski14} \cite[Section~6]{witzel17} that $\phi(\Thomp{G^*})=\phi(G)$, which in turn equals $\phi(G^n)$ for all $n$, so Equation~\eqref{eq:fin} holds.

The matrix group examples are an interesting situation. Theorem~8.1 of \cite{witzel17} says that $\phi(\Thomp{B_*(R)}) \ge \liminf_{n\to\infty}\phi(B_n(R))$ always holds, for any $R$. When $R=\calO_S$, the ring of $S$-integers of a global function field, then $\phi(B_n(\calO_S))=|S|-1$ for $n\ge2$ \cite{bux04}, and Theorem~8.1 of \cite{witzel17} says that $\phi(\Thomp{B_*(\calO_S)})=|S|-1$ as well, so Equation~\eqref{eq:fin} holds as an equality.

Other matrix groups also yield interesting behavior. The \emph{Abels groups} are certain subgroups $\abels_n(\Z[1/p])\le B_n(\Z[1/p])$ that satisfy $\phi(\abels_n(\Z[1/p]))=n-2$. In particular, each group individually is not of type $\F_\infty$, but as $n$ tends to $\infty$, their finiteness lengths do tend to $\infty$, so the right hand side of Equation~\eqref{eq:fin} is actually $\infty$. Theorem~8.10 of \cite{witzel17} says that one can find a cloning system on $(\abels_n(\Z[1/p]))$ and the resulting Thompson-like group $\Thomp{\abels_*(\Z[1/p])}$ is of type $\F_\infty$; hence Equation~\eqref{eq:fin} holds.

When $G_n=S_n^\pm$, as remarked earlier $\Thomp{S_*^\pm}=V^\pm$ is actually isomorphic to $V$, so is of type $\F_\infty$ and Equation~\eqref{eq:fin} holds; this could also be deduced directly from the Stein--Farley complex with an argument very similar to that for $V$. As for $\Vbr^{twist}$, we conjecture that it is of type $\F_\infty$ like all the $B_n^{twist}$, but will leave this for future work.

As a closing remark, it is an interesting question whether there are some nice, general conditions one can impose on a cloning system to ensure that Equation~\eqref{eq:fin} holds. We collect the known results in Table~\ref{tab:fin}, where we see in all cases that Equation~\eqref{eq:fin} holds.

\begin{table}[htb]\label{tab:fin}
\begin{tabular}{c|c|c|c}
  $G_n$ & $\phi(G_n)$ & $\phi(\Thomp{G_*})$ & $\Thomp{G_*}$ \\
  \hline
  ~&~&~&~\\
  $\{1\}$ & $\infty$ & $\infty$ & $F$ \\
  $S_n$ & $\infty$ & $\infty$ & $V$ \\
  $G^n$ & $\phi(G)$ & $\phi(G)$ & $\Thomp{G^*}$ \\
  $B_n$ & $\infty$ & $\infty$ & $\Vbr$ \\
  $PB_n$ & $\infty$ & $\infty$ & $\Fbr$ \\
  $B_n(\calO_S)$ & $|S|-1$ ($n\ge2$) & $|S|-1$ & $\Thomp{B_*(\calO_S)}$ \\
  $\abels_n(\Z[1/p])$ & $n-2$ & $\infty$ & $\Thomp{\abels_*(\Z[1/p])}$ \\
  $S_n^\pm$ & $\infty$ & $\infty$ & $V^\pm$
\end{tabular}
\caption{Comparing finiteness lengths.}
\end{table}

\bibliographystyle{alpha}
\bibliography{cloning_systems_expo}

\end{document}